\documentclass[12pt]{article}
\usepackage{amssymb,amsthm,amsmath,amsfonts,latexsym,tikz,hyperref}
\usepackage[hmargin=1in,vmargin=1in]{geometry}

\newtheorem{thm}{Theorem}[section]
\newtheorem{prop}[thm]{Proposition}
\newtheorem{cor}[thm]{Corollary}
\newtheorem{lem}[thm]{Lemma}
\newtheorem{conj}[thm]{Conjecture}
\newtheorem{exa}[thm]{Example}

\DeclareMathOperator{\tdeg}{totdeg}
\DeclareMathOperator{\sdeg}{sdeg}
\DeclareMathOperator{\Cat}{Cat}
\DeclareMathOperator{\Nar}{Nar}

\newcommand{\ben}{\begin{enumerate}}
\newcommand{\een}{\end{enumerate}}
\newcommand{\ble}{\begin{lem}}
\newcommand{\ele}{\end{lem}}
\newcommand{\bth}{\begin{thm}}
\renewcommand{\eth}{\end{thm}}
\newcommand{\bpr}{\begin{prop}}
\newcommand{\epr}{\end{prop}}
\newcommand{\bco}{\begin{cor}}
\newcommand{\eco}{\end{cor}}
\newcommand{\bcon}{\begin{conj}}
\newcommand{\econ}{\end{conj}}
\newcommand{\bde}{\begin{defn}}
\newcommand{\ede}{\end{defn}}
\newcommand{\bex}{\begin{exa}}
\newcommand{\eex}{\end{exa}}
\newcommand{\barr}{\begin{array}}
\newcommand{\earr}{\end{array}}
\newcommand{\btab}{\begin{tabular}}
\newcommand{\etab}{\end{tabular}}
\newcommand{\beq}{\begin{equation}}
\newcommand{\eeq}{\end{equation}}
\newcommand{\bea}{\begin{eqnarray*}}
\newcommand{\eea}{\end{eqnarray*}}
\newcommand{\bal}{\begin{align*}}
\newcommand{\bce}{\begin{center}}
\newcommand{\ece}{\end{center}}
\newcommand{\bpi}{\begin{picture}}
\newcommand{\epi}{\end{picture}}
\newcommand{\bpp}{\begin{picture}}
\newcommand{\epp}{\end{picture}}
\newcommand{\bfi}{\begin{figure} \begin{center}}
\newcommand{\efi}{\end{center} \end{figure}}
\newcommand{\bprf}{\begin{proof}}
\newcommand{\eprf}{\end{proof}\medskip}
\newcommand{\capt}{\caption}
\newcommand{\bsl}{\begin{slide}{}}
\newcommand{\esl}{\end{slide}}
\newcommand{\bfr}{\begin{frame}}
\newcommand{\efr}{\end{frame}}

\newcommand{\hqed}{\hfill \qed}

\newcommand{\eqed}[1]{$\textcolor{white}{\qed}\hfill{\dil#1}\hfill\qed$}
\newcommand{\eqqed}[1]{$\rule{1ex}{0ex}\hfill{\dil#1}\hfill\qed$}

\newcommand{\hso}[1]{\hspace{-1pt}}
\newcommand{\vs}[1]{\vspace{#1}}
\newcommand{\qmq}[1]{\quad\mbox{#1}\quad}

\newcommand{\Cong}{\equiv}

\newcommand{\case}[4]{\left\{\barr{ll}#1&\mbox{#2}\\#3&\mbox{#4}\earr\right.}
\newcommand{\fl}[1]{\lfloor #1 \rfloor}

\newcommand{\flf}[2]{\left\lfloor\frac{#1}{#2}\right\rfloor}

\def\<{\langle}
\def\>{\rangle}

\newcommand{\ree}[1]{(\ref{#1})}

\newcommand{\ra}{\rightarrow}

\newcommand{\ga}{\gamma}

\newcommand{\ep}{\epsilon}

\newcommand{\ze}{\zeta}

\newcommand{\Ga}{\Gamma}

\newcommand{\1}{{\bf 1}}

\newcommand{\bbC}{{\mathbb C}}

\newcommand{\bbN}{{\mathbb N}}

\newcommand{\bbZ}{{\mathbb Z}}

\newcommand{\cP}{{\cal P}}

\newcommand{\cT}{{\cal T}}

\DeclareMathOperator{\mdeg}{mdeg}
\DeclareMathOperator{\Mod}{mod}

\DeclareMathOperator{\wt}{wt}

\newcommand{\dil}{\displaystyle}

\begin{document}
\pagestyle{plain}

\title{Lucas atoms
}
\author{
Bruce E. Sagan\\[-5pt]
\small Department of Mathematics, Michigan State University,\\[-5pt]
\small East Lansing, MI 48824-1027, USA\\
and\\
Jordan Tirrell\\[-5pt]
\small Department of Mathematics and Computer Science, Washington College,\\[-5pt]
\small Chestertown, MD 21620, USA\\
}

\date{\today\\[10pt]
	\begin{flushleft}
	                                       \\[5pt]
	\end{flushleft}}

\maketitle

\begin{abstract}
Given two variables $s$ and $t$, the associated sequence of Lucas polynomials is defined inductively by $\{0\}=0$, $\{1\}=1$, and
$\{n\}=s\{n-1\}+t\{n-2\}$ for $n\ge2$.  An integer (e.g., a Catalan number) defined by an expression of the form
$\prod_i n_i/\prod_j k_j$ has a Lucas analogue obtained by replacing each factor with the corresponding Lucas polynomial.
There has been  interest in deciding when such expressions, which are a priori only rational functions, are actually polynomials in $s,t$.  
The approaches so far have been combinatorial.  We introduce a powerful algebraic method for answering this question by factoring 
$\{n\}=\prod_{d|n} P_d(s,t)$, where we call the polynomials $P_d(s,t)$ Lucas atoms.  
This permits us to show that the Lucas analogues of the Fuss-Catalan and Fuss-Narayana numbers for all irreducible Coxeter groups are polynomials in $s,t$.
Using gamma expansions, a technique which has recently become popular in combinatorics and geometry,  one can show  that the Lucas atoms have a close relationship with cyclotomic polynomials $\Phi_d(q)$.  
Certain results about the $\Phi_d(q)$ can then be lifted to Lucas atoms.
In particular, one can prove analogues of theorems of Gauss and Lucas, deduce reduction formulas, and evaluate the $P_d(s,t)$ at various specific values of the variables.

\end{abstract}

\section{Introduction}
\label{i}

We will denote the nonnegative integers by $\bbN$.
Let $s,t$ be variables.  Inductively define the $n$th {\em Lucas polynomial}, $\{n\}=\{n\}_{s,t}$, by $\{0\}=0$, $\{1\}=1$, and
\beq
\label{lp}
\{n\}=s\{n-1\} + t\{n-2\}
\eeq
for $n\ge2$.  These polynomials were introduced and studied by Lucas in~\cite{luc:tfn1, luc:tfn2, luc:tfn3}.
This sequence has various interesting specializations.  For example, $\{n\}_{1,1}$ is the $n$th Fibonacci number and 
$\{n\}_{2,-1}=n$.  Furthermore, if one considers a third variable $q$, then a simple induction shows that
\beq
\label{lpq}
\{n\}_{1+q,-q} = 1 + q + q^2 +\cdots + q^{n-1}.
\eeq
This summation is usually denoted $[n]_q$ and is important both in the theory of hypergeometric series and in combinatorics.  This equation will permit us to make a connection between the Lucas sequence and cyclotomic polynomials.

There has been recent interest in studying Lucas analogues of combinatorial constants.  These are are connected via~\ree{lpq} with the well-studied $q$-analogues of such integers.  Suppose we are given an integer defined as a quotient of products $\prod_i n_i/\prod_j k_j$ where the $n_i$ and $k_j$ are positive integers.  The corresponding {\em Lucas analogue} is $\prod_i \{n_i\}/\prod_j \{k_j\}$.  A priori, this is just a rational function of $s$ and $t$.  But often it is actually a polynomial in these variables with nonnegative integer coefficients.  Benjamin and Plott~\cite{bp:caf} gave a complicated combinatorial interpretation for the Lucas analogue of the binomials coefficients, called Lucanomials.  Then Sagan and Savage~\cite{ss:cib} came up with a simpler one which, unfortunately, appeared to be rigid in that their ideas could not be extended to related constants such as the Catalan numbers.  Ekhad~\cite{ekh:ssl} found an algebraic argument to show that since the Lucanomials were in $\bbN[s,t]$, so were the Lucas-Catalans.  Bennett, Carrillo, Machacek, and Sagan~\cite{bcms:cil} gave a combinatorial model in the binomial coefficient case which could be extended to the Catalan numbers for all irreducible Coxeter groups, but they were still not able to apply their methods to various other constants.  As yet unpublished work has also been done by the Algebraic Combinatorics Seminar at the Fields Institute~\cite{abcdl}, Gleb Nenashev~\cite{nen}, and Rao and Suk~\cite{rs:dsp}.

\begin{table}
$$
\barr{c|l|l}
n 	&\{n\} 					&P_n(s,t)\\
\hline
1 	&1 						&1\\
2 	&s 						&s\\
3 	&s^2 + t					&s^2 + t	\\
4 	&s^3 + 2st 				&s^2 + 2t\\
5 	&s^4 + 3 s^2 t + t^2 		&s^4 + 3 s^2 t + t^2 \\
6 	&s^5 + 4 s^3 t + 3 s t^2  	&s^2 + 3t
\earr
$$
\capt{The Lucas polynomials and Lucas atoms for $n\le 6$ \label{lpla}}
\end{table}

We introduce a new and powerful method for proving that Lucas analogues are polynomials with nonnegative integer coefficients. In particular, we will define a new sequence of polynomials $P_n(s,t)$ which will be called Lucas atoms and satisfy
\beq
\label{prod}
\{n\} = \prod_{d|n} P_d(s,t).
\eeq
The first few Lucas polynomials and Lucas atoms are given in Table~\ref{lpla}.  
Given a product of Lucas polynomials $\prod_i \{n_i\}$ its associated {\em atomic decomposition} is the product of Lucas atoms obtained by replacing each $\{n_i\}$ by the corresponding product using~\ree{prod}. 
One of our principal results shows that atomic decompositions function like prime decompositions of integers.  Note that we do not have to consider $P_1(s,t)$ since it is the polynomial $1$.
\bth
\label{factor}
Suppose $f(s,t)=\prod_i \{n_i\}$ and $g(s,t)=\prod_j \{k_j\}$ for certain $n_i,k_j\in\bbN$, and write their atomic decompositions as
$$
f(s,t)=\prod_{d\ge2} P_d(s,t)^{a_d} \qmq{and} g(s,t)=\prod_{d\ge2} P_d(s,t)^{b_d}
$$
for certain powers $a_d, b_d\in \bbN$.
Then $f(s,t)/g(s,t)$ is a polynomial if and only if $a_d\ge b_d$ for all $d\ge2$.  Furthermore, in this case $f(s,t)/g(s,t)$ has nonnegative integer coefficients.
\eth
This result is striking for several reasons.  First of all,  it gives a condition for polynomiality which is not only sufficient but also necessary.   It is also notable that such polynomials must always be in $\bbN[s,t]$.  Thus it is impossible for one of these polynomials to have a coefficient which is $1/2$ or $-3$.  

In the next section, the Lucas atoms are defined and the previous theorem is proved using a connection with cyclotomic polynomials, 
$\Phi_n(q)$.  This correspondence is made through the use of gamma expansions.  These expressions are important in geometry because of a conjecture of Gal and in combinatorics because of their usefulness in proving unimodality results.  See the recent survey of Athanasiadis~\cite{ath:gcg} for more details.
In particular,  $P_n(s,t)$ turns out to be the image of $\Phi_n(q)$ under a map which uses the gamma expansion of the latter.  
It follows that the coefficients of $P_n(s,t)$ are just the absolute values of the gamma coefficients of $\Phi_n(q)$,
In Section~\ref{la} we use Theorem~\ref{factor} to prove that a host of Lucas analogues are in $\bbN[s,t]$, including the Fuss-Catalan and Fuss-Narayana numbers for an arbitrary irreducible Coxeter group.
It is also natural to ask which theorems about the cyclotomic polynomials have counterparts for the Lucas atoms.  Section~\ref{tst} is devoted to showing that theorems of Gauss and Lucas expressing $\Phi_n(q)$ in terms of two squares can be lifted to the Lucas realm.
In Section~\ref{rfsec} we prove reduction formulas for Lucas atoms which reduce their computation to knowing $P_p(s,t)$ for a prime $p$.
Section~\ref{e} contains various evaluations of $P_n(s,t)$ for specific values of $s$ and $t$.  We end with a section of comments and open questions.

\section{Defining Lucas atoms}
\label{dla}

One could define the Lucas atoms $P_n(s,t)$ inductively using~\ree{prod}.  But it will be more useful to obtain them from cyclotomic polynomials.  First, however, we need some definitions about gamma expansions.

Let $p(q)=\sum_{i\ge 0} a_i q^i$ be a nonzero polynomial in $q$ with  coefficients in $\bbC$, the complex numbers.  As usual, the {\em degree} of $p(q)$, $\deg p(q)$, is the largest index $i$ with $a_i\neq 0$.  We will also need the {\em minimum degree} 
$$
\mdeg p(q) = \min\{i \mid a_i\neq 0\}
$$
and {\em total degree}
$$
\tdeg p(q) = \deg p(q) + \mdeg p(q).
$$
For example $p(q) = 2 q + 5 q^2 + 5 q^3 + 2 q^4$ has $\tdeg p(q) = 4 + 1 = 5$.  If $\tdeg p(q) = d$ then we call $p(q)$ 
{\em palindromic} ({\em symmetric} is also used) if $a_i = a_{d-i}$ for all $0\le i\le d$.  It is easy to see that this is equivalent to the equality
\beq
\label{pal}
q^d p(1/q) = p(q).
\eeq
In this case we call $d/2$ the {\em center of symmetry} of $p(q)$.  Our example polynomial is palindromic with center of symmetry $5/2$.  
A straight-forward computation shows that the product of palindromic polynomials is palindromic.  The same is true of linear combinations of palindromic polynomials with the same center of symmetry, but not in general.  

We will need the vector space
$$
\cP_d(q) = \{ p(q) \in \bbC[q] \mid \text{$p(q)$ is palindromic with $\tdeg p(q)=d$}\}\cup\{0\}.
$$
The polynomials
\beq
\label{basis}
(1+q)^d,\ q(1+q)^{d-2},\ q^2(1+q)^{d-4},\ \dots
\eeq
form a basis for $\cP_d$ since they all have different degrees and their leading coefficients equal one.  So if $p(q)\in\cP_d$ then it has 
{\em gamma expansion}
\beq
\label{GamExp}
p(q) = \sum_{j\ge0} \ga_j  q^j (1+q)^{d-2j} 
\eeq
where the scalars $\ga_0,\ga_1,\ga_2,\ldots$ are called the {\em gamma coefficients} of $p(q)$.  Returning to our example
$$
2 q + 5 q^2 + 5 q^3 + 2 q^4 = 0(1+q)^5 + 2q(1+q)^3 -q^2(1+q)
$$
so its gamma coefficients are $0,2,-1$.

To make the connection with the Lucas sequence, an easy inductive proof shows that
\beq
\label{LucExp}
\{n\} = \sum_{j\ge0} a_j s^{n-2j-1} t^j.
\eeq
for certain $a_j\in\bbN$.  Comparison of this expansion with~\ree{GamExp} motivates the following definition.  Consider 
$$
\cP(q) = \bigcup_{d\ge0} \cP_d(q).
$$
Note that the union is disjoint except for the presence of the zero polynomial in all $\cP_d$.  
Define the {\em Gamma map} $\Ga:\cP(q)\ra\bbC[s,t]$ by taking $p(q)$ of the form~\ree{GamExp} to
\beq
\label{GaDef}
\Ga(p(q))=  \sum_{j\ge0} \ga_j  s^{d-2j} (-t)^j.
\eeq
In the next proposition we collect some of the basic properties of this function.
\bpr
\label{GaPr}
The map $\Ga:\cP(q)\ra\bbC[s,t]$ has the following properties
\ben
\item[(a)] If $p(q),r(q)\in\cP(q)$ then
$$
\Ga(p(q) r(q) )=\Ga(p(q))\Ga(r(q)).
$$
\item[(b)]  For any $d$, the restriction of $\Ga$ to $\cP_d(q)$ is linear.
\item[(c)]  The map $\Ga$ is injective.
\item[(d)]  If $\Ga(p(q))=f(s,t)$ then $f(1+q,-q) = p(q)$.
\item[(e)]  If $p(q)\in\bbZ[q]$ then $\Ga(p(q))\in\bbZ[s,t]$.
\een
\epr
\bprf
Parts (a) and (b) follow quickly from the remarks after equation~\ree{pal}.  Also (c) follows from (d) which defines the inverse map on the image of $\Ga$.  For (d) we have, from~\ree{GamExp} and~\ref{GaDef},
$$
f(1+q,-q) = \sum_{j\ge 0} \ga_j s^{d-2j} (-t)^j |_{s=1+q,t=-q} = \sum_{j\ge 0} \ga_j (1+q)^{d-2j} q^j = p(q).
$$
To obtain (e), note that the polynomials in~\ree{basis} are all monic.  So if $p(q)\in\bbZ[q]$, then its gamma coefficients are all integers.  The desired conclusion now follows from the definition of $\Ga$.
\eprf

To define the Lucas atoms, we first recall some simple facts about the cyclotomic polynomials.  The {\em $n$th cyclotomic polynomial} is 
$$
\Phi_n(q) = \prod_\ze (q-\ze)
$$
where the product is over all primitive $n$th roots of unity.  Since $\ze$ is a primitive $n$th root if and only if $1/\ze$ is, and the constant coefficient of $\Phi_n(q)$ is one for $n\ge2$, it follows from equation~\ree{pal} that $\Phi_n(q)$ is palindromic for that range of $n$.  So define the {\em $n$th Lucas atom} as $P_1(s,t)=1$ and 
$$
P_n(s,t) =\Ga(\Phi_n(q))
$$
for $n\ge2$.  The basic properties of $P_n(s,t)$ are as follows.
\bpr
\label{P_nPr}
For all $n\ge1$ we have
\ben
\item[(a)] $\dil \{n\} = \prod_{d|n} P_d(s,t)$, and 
\item[(b)] $P_n(s,t)\in\bbN[s,t]$.
\een
\epr
\bprf
(a)   It is  well known and easy to prove from the definitions that
\beq
\label{PhiProd}
q^n-1 =\prod_{d|n} \Phi_d(q).
\eeq
It follows that
$$
1+q+q^2+\cdots + q^{n-1} = \prod_{d|n\atop d\ge2} \Phi_d(q).
$$
So applying Proposition~\ref{GaPr} (a) and using the fact that $P_1(s,t)=1$ we have
$$
\Ga(1+q+q^2+\cdots + q^{n-1}) =\prod_{d|n} P_d(s,t).
$$
But from equation~\ree{lpq} as well as Proposition~\ref{GaPr} (c) and (d) we have that the left side of the previous equation is $\{n\}$.

\medskip

(b)  Since the leading coefficient of $\{n\}$ is one, an easy induction using part (a) shows that the same is true of the $P_n(s,t)$.  A second induction based on (a) now shows that all the coefficients of $P_n(s,t)$ are integers.   For nonnegativity, it suffices to show that, for  $n\ge3$, the polynomial $P_n(s,t)$ can be written as a product of factors of the form $s^2+at$ where $a>0$.  (Nonnegativity for $n\le2$ is clear.)  Consider any root $\ze$ of $\Phi_n(q)$.  Then the complex conjugate $\bar{\ze}$ is also a root, and $\Phi_n(q)$ has a factor
$$
(q-\ze)(q-\bar{\ze}) = q^2 -2b q + 1 = (q+1)^2 - (2b+2) q
$$
where $b$ is the real part of $\ze$.  Since $n\ge3$ we have $a:=2b+2>0$.  Using Proposition~\ref{GaPr} (a) shows that
$$
\Ga((q+1)^2 - (2b+2) q) = s^2+at
$$
is a factor of $P_n(s,t)$ as desired.
\eprf

We now have all the tools necessary to prove Theorem~\ref{factor}.

{\em Proof (of Theorem~\ref{factor}).}
Clearly if $a_d\ge b_d$ for all $d\ge2$ then $f(s,t)/g(s,t)$ is a polynomial.  And since the $P_d(s,t)$ all have nonnegative integer coefficients from the previous proposition, the inequalities show that the same is true of the quotient since it is a product of atoms.  So it remains to show that $f(s,t)/g(s,t)$ being a polynomial implies that $a_d\ge b_d$ for $d\ge2$.  It is clear that this holds for $d=2$ since
$P_2(s,t)=s$.  So suppose $d\ge3$.   Now $f(s,1)/g(s,1)$ is a polynomial in $s$.  And from the proof of Proposition~\ref{P_nPr} (b), we see that the roots of $P_d(s,1)$ are all of the form $\pm\sqrt{-2b-2}$ where $b$ is the real part of a primitive $d$th root of unity.  It follows that no two of these polynomials in $s$ have a common root.  So the polynomialty of $f(s,1)/g(s,1)$  implies  $a_d\ge b_d$ for all $d\ge2$.
\hqed

\section{Lucas analogues}
\label{la}

We will now use Theorem~\ref{factor} to show that a large number of Lucas analogues are polynomials with nonnegative integer coefficients.  We will start with the binomial coefficients, then consider various types of Fuss-Catalan numbers including those associated with irreducible Coxeter groups, and finally look at Fuss-Narayana numbers.

We first need to consider the Lucas factorization of the {\em Lucatorial}
$$
\{n\}! = \{1\} \{2\}\cdots\{n\}.
$$
To describe the factorization we will need the {\em floor} or {\em round-down function} $\fl{q}$ which is the largest integer less than or equal to the rational number $q$.  Given a product $f(s,t)$ of Lucas polynomials, let
$$
\log_d f(s,t) = \text{the power of $P_d(s,t)$ in its Lucas factorization.}
$$
The subscript  will be omitted if $d$ is clear from context or is generic and fixed.
\ble
\label{flf}
For $d\ge2$ we have
$$
\log_d \{n\}! = \fl{n/d}.
$$
Furthermore, for integers $m,n,d$
$$
\flf{m}{d} + \flf{n}{d} \le \flf{m+n}{d}.
$$
\ele
\bprf
We have that $P_d$ is a factor of $\{j\}$ if and only if $d|j$.  So the number of such factors in $\{n\}!$ is $\fl{n/d}$.  The inequality is well known so we omit the proof.
\eprf

Now for $0\le k\le n$, we define the {\em Lucanomial}
$$
\left\{n\atop k\right\} = \frac{\{n\}!}{\{k\}!\{n-k\}!}.
$$
\bth
\label{lucanom}
For $0\le k\le n$ we have $\left\{n\atop k\right\}\in\bbN[s,t]$.
\eth
\bprf
Applying the previous lemma gives, for any $d\ge 2$,
$$
\log_d (\{k\}! \{n-k\}!) = \flf{k}{d}+\flf{n-k}{d}\le \flf{n}{d}= \log_d\{n\}!.
$$
So we are done by Theorem~\ref{factor}.
\eprf

We will now consider various types of Catalan numbers.  Given positive integers $a,b$ with $\gcd(a,b)=1$ the corresponding 
{\em rational Catalan number} is
$$
\Cat(a,b) = \frac{1}{a+b} \binom{a+b}{a}.
$$
One obtains the usual Catalan numbers by letting $a=n$ and $b=n+1$.  The corresponding Lucas analogue is
$$
\Cat\{a,b\}= \frac{1}{\{a+b\}} \left\{a+b \atop a\right\}.
$$
The Algebraic Combinatorics Seminar at the Fields Institute~\cite{abcdl} was the first to prove that the $q$-analogue of $\Cat(a,b)$ is a polynomial in $q$ and their method works as well for the Lucas analogue.  This proof is also algebraic and is presented in~\cite{bcms:cil}.  A combinatorial proof has yet to be found.
\bth
\label{RatCat}
If $\gcd(a,b)=1$ then $\Cat\{a,b\}\in\bbN[s,t]$.
\eth
\bprf
There are two cases.  If $d$ does not divide $a+b$ then $\log\{a+b\}=0$ and so the result follows from the previous theorem.
If $d$ divides $a+b$ then $d$ divides at most one of $a$ and $b$ since $\gcd(a,b)=1$.  It follows that either 
$\fl{a/d}=\fl{(a-1)/d}$ or $\fl{b/d}=\fl{(b-1)/d}$.  So for either possibility
$$
\log(\{a+b\}\{a\}!\{b\}!) \le 1 + \flf{a+b-1}{d} = \flf{a+b}{d} = \log\{a+b\}!
$$
by Lemma~\ref{flf}.  Theorem~\ref{factor} completes the proof.
\eprf

\begin{table}
$$
\barr{l|l}
W 	& d_1,\dots,d_n\\
\hline
A_n 	& 2,3,4,\dots,n+1\\
B_n & 2,4,6,\dots,2n\\
D_n & 2,4,6,\dots,2n-2,n\\
E_6 & 2,5,6,8,9,12\\
E_7 & 2,6,8,10,12,14,18\\
E_8 & 2,8,12,14,18, 20, 24, 30\\
F_4 & 2, 6, 8, 12\\
H_3 & 2, 6, 10\\
H_4 & 2, 12, 20, 30\\
I_2(m) & 2, m
\earr
$$
\capt{The finite irreducible Coxeter groups and their degrees~\label{CoxTab}}
\end{table}

Let $W$ be a finite irreducible Coxeter group with degrees $d_1<\dots<d_n$.  A list of these groups and their degrees is given in Table~\ref{CoxTab}.  If $k$ is a positive integer then $W$ has corresponding {\em Fuss-Catalan number}
$$
\Cat^{(k)} W = \prod_{j=1}^n \frac{d_j+ k d_n}{d_j}.
$$
The study of these constants and related ideas has come to be known as ``Coxeter-Catalan combinatorics."  See the memoir of Armstrong~\cite{arm:gnp} for more information.  The corresponding Lucas analogue is
$$
\Cat^{(k)} \{W\} = \prod_{j=1}^n \frac{\{d_i + k d_n\}}{\{d_n\}}.
$$
When referring to a specific $W$, we put the curly brackets around the subscript giving the rank, e.g., $\Cat^{(k)} B_{\{n\}}$.
\bth
\label{CatThm}
For all finite irreducible Coxeter groups $W$ and all  positive integers $k$ we have $\Cat^{(k)} \{W\}\in\bbN[s,t]$.
\eth
\bprf
We note that for the classical types $A_n, B_n, D_n$ Bennett et al.~\cite{bcms:cil} were able to prove this result by combinatorial arguments.  It remains open to do the same for the exceptional groups.  We will proceed group by group.

\medskip

{\bf Type $A_{n-1}$.} In this case we can express the Fuss-Catalan analogue in terms of the rational Catalan analogue since
$$
\Cat^{(k)} A_{\{n-1\}} =\frac{\{(k+1)n\}!}{\{n\}!\{kn+1\}!}=\Cat\{n,kn+1\}.
$$
So the result follows from the previous theorem.

\medskip

{\bf Type $B_n$.}
In type $B_n$  one can cancel powers of two from the numerator and denominator and so express $\Cat^{(k)} B_n$ as a binomial coefficient.  But one can no longer do this when each factor is replaced by the corresponding Lucas polynomial.  Instead we will consider a generalization of $\Cat^{(k)} B_{\{n\}}$.  Given an integer $m\ge1$ we let
$$
\{n:m\}! = \{m\} \{2m\} \cdots \{nm\}
$$
and define an {\em $m$-divisible Lucanomial} to be
\beq
\label{n:m}
{n:m \brace k:m} =\frac{\{n:m\}!}{\{k:m\}!\{n-k:m\}!}.
\eeq
So we have the special case
$$
\Cat^{(k)} B_{\{n\}} = {(k+1)n: 2 \brace n:2}.
$$
To show that~\ree{n:m}  is in $\bbN[s,t]$, note that $P_d$ divides terms at interval $d/\gcd(d,m)$ in 
$\{n:m\}!$.  The rest of the proof is much the same as for the Lucanomials and so is omitted.

\medskip

{\bf Type $D_n$.}
We have
\begin{align*}
\Cat^{(k)} D_{\{n\}} &= \frac{\{n+2(n-1)k\}}{\{n\}}{(k+1)(n-1):2 \brace n-1:2}\\[20pt]
&= \frac{\{n+2(n-1)k\}}{\{n\}} \cdot\frac{\{2+2(n-1)k\}\{4+2(n-1)k\}\cdots\ \{2(k+1)(n-1)\}}{\{n-1:2\}!}.
\end{align*}
Given $d\ge2$ there are two cases.  If $d$ does not divide $n$ then the factors of $P_d$ in the denominator all occur inside the $2$-divisible Lucanomial and so cancel out as for type $B_n$.  

Now suppose $d|n$.  In any product of the form $\{2l\}\{2l+2\}\cdots\{2m\}$, the Lucas atom $P_d$ will divide terms at intervals of length $d'=d/\gcd(d,2)$.  Since $d|n$ we have that $P_d$ will appears in exactly $n/d' - 1$ factors in $\{n-1:2\}!$, giving a total of $n/d'$ times in the denominator of the Fuss-Catalan analogue.
If $d$ does not divide $2(n-1)k$ then $P_d$ will divide $n/d'$ terms in the numerator of the last fraction of the above displayed equation and we will be done.  If $d|2(n-1)k$ then $P_d$ will only divide $n/d'-1$ terms in that product, but will also divide $\{n+2(n-1)k\}$ in the numerator, giving the required number of $n/d'$ copies.

\medskip

{\bf Type $I_2(m)$.}
We have
$$
\Cat^{(k)} I_{\{2\}}(m) =\frac{\{km+2\}\{(k+1)m\}}{\{2\}\{m\}}.
$$
If $P_d$ appears as a factor in the denominator for $d\ge3$ then we must have $d|m$.  It follows that $d|(k+1)m$ and so is canceled by the corresponding factor in the numerator.  If $d=2$ then there are two cases.  If $m$ is even then similar consideration show that $P_2^2$ appears in both the denominator and the numerator.  If $m$ is odd then the denominator only has one $P_2$.  In the numerator, that factor will appear in $\{km+2\}$ if $k$ is even or $\{(k+1)m\}$ if $k$ is odd.

\medskip

{\bf The exceptional types.}
For the exceptional $W$, we do not need to consider an infinite number of values of $k$.  This is because whether a given $P_d$ divides a factor $\{a+bk\}$  in the numerator depends only on the congruence class of $k$ modulo $d$.  And the number of choices for $d$ is limited by the factors in the denominator.  But those factors do   not depend on $k$ and so there are only finitely many choices.  In fact, these demonstrations are so straightforward that they can easily be done by hand.  So we will only illustrate the procedure in a particular example.  Consider 
$$
\Cat^{(k)} H_{\{4\}} = \frac{\{2+30k\}\{12+30k\}\{20+30k\}\{30+30k\}}{\{2\}\{12\}\{20\}\{30\}}.
$$
Now $P_4$ is in the factorization of $\{12\}$ and $\{20\}$ in the denominator, so we must show it is also appears in the expansion of two of the factors in the numerator regardless of $k$.  Reducing modulo $4$, we see that it suffices to look at $P_4$ factors of 
$$
\{2+2k\}\cdot\{2k\}\cdot\{2k\}\cdot\{2+2k\}=\{2k\}^2\cdot\{2(k+1)\}^2.
$$
So if $k$ is even, then $P_4^2$ appears in $\{2k\}^2$, while if $k$ is odd then it divides $\{2(k+1)\}^2$. 
\eprf

Let $W$ be a finite irreducible Coxeter group of rank $n$, and let $k,i$ be integers with $k$ positive and $0\le i\le n$.  The corresponding {\em Fuss-Narayana numbers} are $\Nar^{(k)}(W,i)$ which count the number of $k$-multichains in the noncrossing partition poset  of $W$ whose bottom element has rank $i$.  It can be proved that these are always polynomials in $k$~\cite[Theorem 3.5.5]{arm:gnp}.  However, there does not seem to be a simple product formula for them which holds for all $W,k,i$.  However, when $i=1$ we have
\beq
\label{NarPro}
\Nar^{(k)} W := \Nar^{(k)}(W,1) = n \prod_{j=1}^{n-1} \frac{k d_n-d_j+2}{d_j}.
\eeq
Also, for all $i$,
\begin{align*}
 \Nar^{(k)}(A_{n-1},i) &=\frac{1}{n}\binom{n}{i}\binom{kn}{n-i-1},\\[5pt]
 \Nar^{(k)}(B_n,i) &=\binom{n}{i}\binom{kn}{n-i},\\[5pt]
 \Nar^{(k)}(D_n,i) &=\binom{n}{i}\binom{k(n-1)}{n-i}+\binom{n-2}{i}\binom{k(n-1)+1}{n-i}.
\end{align*}
We denote the Lucas analogues by replacing $W$ with $\{W\}$ for a general Coxeter group, and replacing a subscript $n$ by $\{n\}$ in types $A$, $B$, and $D$.  We note that Nenashev~\cite{nen}  has discovered a combinatorial interpretation for the polynomials $\Nar^{(1)}(A_{\{n-1\}},i)$.
\bth
\label{NarThm}
For all finite irreducible Coexter groups and all positive integers $k$ we have $\Nar^{(k)} \{W\} \in\bbN[s,t]$.  This is also true of 
$\Nar^{(k)}(A_{\{n-1\}},i)$,  $\Nar^{(k)}(B_{\{n\}},i)$, and $\Nar^{(k)}(D_{\{n\}},i)$ for all $i$.
\eth
\bprf
The $\Nar^{(k)} \{W\}$ are taken care of in much the same way as the Fuss-Catalan analogues for the exceptional groups.  
And in type $B$ and $D$ the analogues are just sums and products of Lucasnomials.
So we will only give details for
$$
\Nar^{(k)}(A_{\{n-1\}},i)=\frac{1}{\{n\}}\left\{n\atop i\right\}\left\{kn\atop n-i-1\right\}.
$$
If $d$ does not divide $n$ then $P_d$ can only appear in the denominators of the Lucanomials and so the inequality in 
Theorem~\ref{factor} holds for these $d$ by Theorem~\ref{lucanom}.  If $d|n$ then it can not divide both $i$ and $n-i-1$.
The demonstration is now completed as in the proof of Theorem~\ref{RatCat}.
\eprf

\section{Two square theorems}
\label{tst}

In this section we will find Lucas analogues of theorems of Gauss and Lucas which express (appropriately modified) cyclotomic polynomials in terms of squares of two other polynomials.  This will turn out to be easy to do by applying the function $\Ga$ defined by~\ree{GaDef}.

The result of Lucas~\cite[pp.\ 309--315, p.\ 443]{rie:pnc} is as follows.  In it, $\phi(n)$ denotes the Euler totient function.
\bth[Lucas' formula]
\label{LucThm}
If $n\ge 5$ is odd and square-free, then there are polynomials $C_n(q)$ and $D_n(q)$ such that
$$
\Phi_n\left((-1)^{(n-1)/2}q\right) = C_n^2(q) - n q D_n^2(q).
$$
If $n\ge4$ is even and square-free then there are polynomials $C_n(q)$ and $D_n(q)$ such that
$$
\Phi_{2n}(q) = C_n^2(q) - n q D_n^2(q).
$$
In both cases
\ben
\item $C_n(q),D_n(q)\in\bbZ[q]$,
\item $\deg C_n(q) = \phi(n)/2$ and $\deg D_n = \phi(n)/2- 1$,
\item $C_n(q)$ and $D_n(q)$ are both palindromic.
\een
\eth

To state the analogous result for Lucas atoms we define, for $f(s,t)\in\bbC[s,t]$, 
$$
\sdeg f(s,t) = \text{largest power of $s$ in $f(s,t)$.}
$$
\bth
\label{LucAna}
If $n\ge 5$ is square-free and satisfies $n\Cong 1\ (\Mod 4)$, then there are polynomials $G_n(s,t)$ and $H_n(s,t)$ such that
\beq
\label{PGH}
P_n(s,t) = G_n^2(s,t) + n t H_n^2(s,t).
\eeq
If $n\ge4$ is even and square-free, then there are polynomials $G_n(s,t)$ and $H_n(s,t)$ such that
$$
P_{2n}(s,t) = G_n^2(s,t) + n t H_n^2(s,t).
$$
In both cases
\ben
\item $G_n(s,t), H_n(s,t)\in\bbZ[s,t]$,
\item $\sdeg G_n(s,t) = \phi(n)/2$ and $\sdeg H_n(s,t) = \phi(n)/2- 1$.\hqed
\een
\eth
\bprf
We will only prove the statement about odd $n$ as the one for even values is obtained similarly.   Since $n\Cong 1\ (\Mod 4)$ we have from Theorem~\ref{LucThm} that 
\beq
\label{PhiCD}
\Phi_n(q)=\Phi_n\left((-1)^{(n-1)/2}q\right) = C_n^2(q) - n q D_n^2(q).
\eeq
Since $n\ge5$ we know that $\Phi_n(q)$ is palindromic with $\sdeg P_n(s,t)=\deg \Phi_n(q) = \phi(n)$.  From the given facts about $C_n(q)$ and $D_n(q)$ we see that $C_n^2(q)$ and $q D_n^2(q)$ are both palindromic of total degree $\phi(n)$.  So, from Proposition~\ref{GaPr} (a) and (b), we can apply $\Ga$ to both sides of~\ree{PhiCD} and obtain~\ree{PGH}.  The fact that the  polynomials $G$ and $H$ have integer coefficients is a consequence of Proposition~\ref{GaPr} (e).  And the statement about their degrees follows directly from the definition of $\Ga$.
\eprf

One might wonder if it is possible to get an analogue of Lucas' formula when $n$ is square-free and congruent to $3$ modulo $4$.
However, one does not seem to exist.  For example, we have
$$
P_7(s,t) =s^6 + 5 s^4 t + 6 s^2 t^2 + t^3
$$
If the desired $G_7(s,t)$ and $H_7(s,t)$ did exist, then the term $t^3$ in $P_7(s,t)$ could not come from $G_7^2$ because 
of the  odd power of $t$.  But $t^3$ could also not arise from $7 t H_7^2$ since $H_7\in\bbZ[s,t]$ and so every term in the product has coefficient divisible by $7$.

We will now consider the formula of Gauss~\cite[Articles 356--357]{gau:da}.  To state it, we define a polynomial $p(q)=\sum_i a_i q^i$ with $\tdeg p(q)=d$ 
to be {\em anti-palindromic} if $a_i = -a_{d-i}$ for  all $0\le i\le d$.
\bth[Gauss' formula]
\label{GauThm}
If  $n\ge 5$ is odd and square-free, then there are polynomials $A_n(q)$ and $B_n(q)$. such that
$$
4\Phi_n(q) = A_n^2(q) - (-1)^{(n-1)/2} n q^2 B_n^2(q)
$$
where 
\ben
\item $A_n(q),B_n(q)\in\bbZ[q]$,
\item $\deg A_n(q) = \phi(n)/2$ and $\deg B_n = \phi(n)/2 - 2$,
\item If $n\Cong 1\ (\Mod 4)$  then $A_n(q)$ and $B_n(q)$ are  palindromic.
\item If $n\Cong 3\ (\Mod 4)$ then $A_n(q)$ is antipalindromic and $B_n$ is palindromic if $n$ is prime, or vice-versa if $n$ is composite.
\een
\eth

Again, only the case when $n\Cong 1\ (\Mod 4)$ seems to have a Lucas analogue.  The proof of the next result is close enough to that of  Theorem~\ref{LucAna} that we leave it to the reader.
\bth
If  $n\ge 5$ is  square-free and satisfies $n\Cong 1\ (\Mod 4)$, then there are polynomials $E_n(s,t)$ and $F_n(s,t)$. such that
$$
4P_n(s,t) = E_n^2(s,t) - n t^2 F_n^2(s,t)
$$
where 
\ben
\item $E_n(s,t),F_n(s,t)\in\bbZ[s,t]$,
\item $\sdeg E_n(s,t) = \phi(n)/2$ and $\sdeg F_n(s,t)= \phi(n)/2 - 2$.\hqed
\een
\eth

\section{Reduction formulas}
\label{rfsec}

The reduction formulas permit the calculation of $\Phi_n(q)$ in terms of $\Phi_m(q)$ for $m<n$.  And these computations are  done over the integers rather than the complex numbers.  The following reductions are all easy to prove directly from the definition of $\Phi_n(q)$ and properties of primitive roots of unity.
\bth[Reduction formulas]
\label{rf}
Let $n$ be a positive integer and $p$ be a prime not dividing $n$.
\ben
\item[(a)] We have
$$
\Phi_p(q) = [p]_q = 1 + q + q^2 +\cdots + q^{p-1}.
$$
\item[(b)] If $m\ge2$ then 
$$
\Phi_{p^m n}(q) = \Phi_{pn}(q^{p^{m-1}}).
$$
\item[(c)]  For all $p$ we have
$$
\Phi_{pn}(q) = \frac{\Phi_n(q^p)}{\Phi_n(q)}.
$$
And for $p=2$ we also have

\eqqed{
\Phi_{2n}(q) = \Phi_n(-q).
}
\een
\eth

So given any $n$, we can use part (b) to reduce the calculation of $\Phi_n(q)$ to that of the radical (square-free part) of $n$. Then part (c)  turns computation for the radical into knowing $\Phi_p(q)$ for primes $p$.  And for these we have an explicit formula in part (a).

\bfi
\begin{tikzpicture}
\draw (0,0) grid (3,1);
\fill (.5,.5) circle (.1);
\fill (1.5,.5) circle (.1);
\fill (2.5,.5) circle (.1);
\end{tikzpicture}
\qquad
\begin{tikzpicture}
\draw (0,0) grid (3,1);
\draw (.5,.5)--(1.5,.5);
\fill (.5,.5) circle (.1);
\fill (1.5,.5) circle (.1);
\fill (2.5,.5) circle (.1);
\end{tikzpicture}
\qquad
\begin{tikzpicture}
\draw (0,0) grid (3,1);
\draw (1.5,.5)--(2.5,.5);
\fill (.5,.5) circle (.1);
\fill (1.5,.5) circle (.1);
\fill (2.5,.5) circle (.1);
\end{tikzpicture}
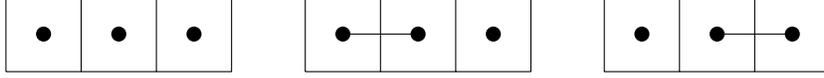
\capt{The tilings in $\cT(3)$}
\label{cT(3)}
\efi

It does not seems as if one can find analogues for these formulas merely by applying  $\Ga$.  The problem is that the necessary substitutions do not seem to behave well with respect to this map.  Instead, we will need a number of lemmas.  For some of them, it will be convenient to use a combinatorial description of $\{n\}$ in terms of tilings.  
For more information about this approach, see the book of Benjamin and Quinn~\cite{bq:prc}.
Consider a row of $n$ boxes.  A {\em tiling}, $T$, of  this row is a covering of the boxes with disjoint tiles where each tile covers two boxes (called a {\em domino}) or one box (called a {\em monomino}). Let $\cT(n)$ denote the set of such tilings.  The set $\cT(3)$ is displayed in Figure~\ref{cT(3)}.  Give a single tiling $T$ the weight
$$
\wt T = s^{\text{number of monominos in $T$}} t^{\text{number of dominos in $T$}}.
$$
Also weight any set $\cT$ of tilings by
$$
\wt\cT=\sum_{T\in\cT} \wt T.
$$
Returning to Figure~\ref{cT(3)} we see that $\wt(\cT(3))=s^3+2st=\{4\}$.  This illustrates a general result which is easy to prove by induction and gives a combinatorial explanation for equation~\ree{LucExp}.
\ble
\label{cT(n)}
For all $n\ge 1$ we have

\eqqed{
\{n\} = \wt(\cT(n-1)).
}
\ele

From this result, we get our Lucas analogue of Theorem~\ref{rf} (a).
If $f$ is a polynomial in $s,t$ then let $[s^i t^j]f$ be the coefficient of $s^i t^j$ in $f$.
\bco
\label{P_p}
For $n\ge1$ we have
$$
\{n\}=\sum_{k\ge 0} \binom{n-k-1}{k} s^{n-2k-1}t^k.
$$
So if $p$ is prime then
$$
P_p(s,t) = \sum_{k\ge 0} \binom{p-k-1}{k} s^{p-2k-1}t^k.
$$
\eco
\bprf
The second statement follows from the first and the fact that for a prime $p$ we have $\{p\}=P_1 P_p = P_p$.
To prove the first, from the previous lemma, $[s^{n-2k-1} t^k]\{n\}$ is the number of tilings of $\cT(n-1)$ with $k$ dominoes and $n-2k-1$ monominoes.  But the number of ways to do this is the number of ways of choosing $k$ dominoes from a total of $n-k-1$ tiles, giving the desired binomial coefficient.
\eprf

The odd primes and $2$ will take different roles in our investigation.  So we will need the following result.

\ble
\label{P_2pThm}
If $p$ is an odd prime then
\beq
\label{P_2p}
P_{2p}(s,t) = \sum_{k\ge 0} \left[\binom{p-k}{k}+\binom{p-k-1}{k-1}  \right]s^{p-2k-1} t^k,
\eeq
and
\beq
\label{PvsL}
s P_{2p}(s,t) = \{p+1\}+t\{p-1\}.
\eeq
\ele
\bprf
The second equation follows from the first the previous corollary.  We now prove the first.
By Proposition~\ref{P_nPr} (a), it suffices to  let $Q_{2p}$ be the right-hand side of~\ree{P_2p} and show that
$$
\{2p\} = P_1 P_2 P_p Q_{2p}
$$
But, again using the previous corollary,
$$
 P_1 P_2 P_p Q_{2p}
= 
s\left(\sum_{i\ge0} \binom{p-i-1}{i} s^{p-2i-1} t^i \right)
\left(\sum_{j\ge0}  \left[\binom{p-j}{j}+\binom{p-j-1}{j-1}  \right]s^{p-2j-1} t^j\right)
$$
So
\beq
\label{R2p}
[s^{2p-2k-1} t^k ]P_1 P_2 P_p Q_{2p}
= \sum_{i+j=k} \binom{p-i-1}{i}\binom{p-j}{j}+  \sum_{i+j=k-1}\binom{p-i-1}{i}\binom{p-j-2}{j}.
\eeq
Using Corollary~\ref{P_p} yet again
$$
[s^{2p-2k-1} t^k ]\{2p\}= \binom{2p-k-1}{k}.
$$

To show equality of the right-hand sides of the previous two equations note that  the single binomial coefficient is the number of tilings  of  $2p-1$  squares with $k$ dominoes.  These tilings are of two types: those with no domino between the $(p-1)$st and $p$th squares and those where these two squares contain a domino.  The first sum in~\ree{R2p} counts the first set of tilings because they can be formed by concatenating a tiling of $p-1$ squares having $i$ dominoes with a tiling of $p$ squares have $j$ dominoes where $i+j=k$.  Similarly, the second sum enumerates the second set of tilings since after the given domino is removed then one is left with a tiling of $p-2$ squares and a tiling of $p-1$ squares with a total of $k-1$ dominoes.
\eprf

Our goal now is to prove an analogue of Theorem~\ref{rf} (c) for Lucas atoms.  We still need several lemmas.  The next result is simple to prove using an argument like that in the last paragraph of the previous demonstration.  So we omit the proof.
\ble
\label{m+n}
For $m,n\ge0$  we have

\eqqed{
\{m+n\} = \{m+1\}\{n\}+t\{m\}\{n-1\}.
}
\ele

We use the notation $M$ or$D$ for a monomino or domino tile, respectively.  Also, $ST$ will denote the concatenation of tilings $S$ amd $T$ and we will use multiplicity notation such as $T^2$ for the concatenation of $T$ with itself.
We also let $\{n\}=0$ if $n\le 0$.
Define the {\em sign} of an integer $m$ to be
$$
\ep(m) = \case{-1}{if $m$ is even,}{+1}{if $m$ is odd.}
$$

\ble
\label{TileSq}
For $m\ge1$ we have
$$
\{m\}^2=\{m-1\}\{m+1\}+\ep(m) t^{m-1}
$$
\ele
\bprf
We will give a proof when $m$ is odd as the other case is similar.  Let $n=m-1$.  It suffices to find a weight-preserving bijection
$$
f:[\cT(n)\times \cT(n)]^- \ra \cT(n-1)\times\cT(n+1)
$$
where $[\cT(n)\times \cT(n)]^-$ is $\cT(n)\times \cT(n)$ with the pair $(D^{n/2}, D^{n/2})$ removed.
Label the $n$ squares left to right from $1$ to $n$.
Given a pair $(S,T)$ in the domain, consider the largest index $i\ge0$ such that only dominoes cover squares of index less than or equal to 
$i$ in both $S$ and $T$.  So $i$ is even and write $S=D^{i/2} S'$ and $T=D^{i/2} T'$.  Since $(S,T)\neq (D^{n/2}, D^{n/2})$ the tilings $S',T'$ are nonempty.  If $S'=MS''$ for some $S''$ then let $f(S,T) = (D^{i/2}S'',D^{i/2}MT')$.  If $S'=DS''$ then, by maximality of $i$, we must have $T'=MT''$ for some $T''$.  In this case let $f(S,T) = (D^{i/2}MS'',D^{i/2}DT'')$.  Clearly this map preserves weight.  And its inverse is easy to construct, so it is bijective.
\eprf

The next lemma can be thought of as a combination of the previous two.

\ble
\label{mnrr}
If $n\ge 2m$ then
$$
\{n\}=\left(\{m+1\}+t\{m-1\}\right)\{n-m\}+\ep(m) t^m\{n-2m\}.
$$
\ele
\bprf
We induct on $n$, assuming $m$ is odd as the even case is similar.  For the base cases, first consider $n=2m$.  So we wish to prove
$$
\{2m\}= \{m+1\}\{m\} + t\{m-1\}\{m\}
$$
which follows by letting $m=n$ in Lemma.~\ref{m+n}.  For the other base case, suppose $n=2m+1$ and compute the right-hand side of the equality using Lemma~\ref{TileSq} and then Lemma~\ref{m+n}
$$
\{m+1\}^2 + t\{m-1\}\{m+1\} +t^m
=\{m+1\}^2 + t(\{m\}^2-t^{m-1})+t^m
=\{2m+1\}.
$$

For the induction step, we use the defining recursion for the Lucas sequence several times on the right-hand side of the desired equation, letting $A=\{m+1\}+t\{m-1\}$ for readability,
$$
\barr{l}
A\{n-m\}+t^m\{n-2m\}\\[5pt]
\qquad=A(\{n-m-1\}+t\{n-m-2\})+t^m(\{n-2m-1\}+t\{n-2m-2\})\\[5pt]
\qquad=(A\{n-m-1\} + t^m\{n-2m-1\}) + t (A\{n-m-2\}+t^m\{n-2m-2\})\\[5pt]
\qquad= \{n-1\}+t\{n-2\}\\[5pt]
\qquad= \{n\}
\earr
$$
which is what we wished to show.
\eprf

We have one last identity to prove before demonstrating our first main theorem of this section.  Note that we can unify the two cases in the following results by using the fact that for $p$ prime we have, by equation~\ree{PvsL},
\beq
\label{ppm1}
\{p+1\}+t\{p-1\}
=\case{s^2+2t}{if $p=2$,}{sP_{2p}}{if $p\ge3$.\rule{0pt}{15pt}}
\eeq
But because of the subscripts, it is easier to read these results in the format we present.
\ble
\label{pn:decomp}
If $p$ is prime then for all $n\ge0$ we have
$$
\{pn\}=\case{\{p\}\cdot \{n\}_{s^2+2t, -t^2}}{if $p=2$,}{\{p\}\cdot \{n\}_{sP_{2p}, t^p}}{if $p\ge3$.\rule{0pt}{15pt}}
$$
\ele
\bprf
We will do the case for odd primes as $p=2$ is similar.
Induct on $n$.  The identity is easy to check when $n=0,1$.  For $n\ge2$ we use in turn the recursion defining the Lucas sequence, induction, equation~\ree{ppm1}, and Lemma~\ref{mnrr} (with $n$ replaced by $pn$ and $m$ replaced by $p$) to obtain
\begin{align*}
\{p\}\cdot \{n\}_{s P_{2p}, t^p}
&=\{p\}\left( sP_{2p}\cdot \{n-1\}_{s P_{2p}, t^p}  + t^p \cdot \{n-2\}_{s P_{2p}, t^p}\right)\\[5pt]
&=sP_{2p} \cdot \{pn-p\} + t^p \cdot\{pn-2p\}\\[5pt]
&=(\{p+1\}+t\{p-1\})\cdot \{pn-p\} + t^p \cdot\{pn-2p\}\\[5pt]
&=\{pn\}
\end{align*}
as desired.
\eprf

We can finally prove our analogue of Theorem~\ref{rf} (c).  
\bth
\label{pn}
If $n\ge2$ is a positive integer and $p$ is a prime not dividing $n$, then
$$
P_{pn}(s,t) = \case{ \frac{\dil P_n(s^2+2t,-t^2)}{\dil P_n(s,t)}}{if $p=2$,}{\frac{\dil P_n(s P_{2p}, t^p)}{\dil P_n(s,t)}}{if $p\ge3$.\rule{0pt}{25pt}}
$$
\eth
\bprf
We assume $p$ is odd as $p=2$ is similar.   We also continue to use $P_n$ as an abbreviation for $P_n(s,t)$, but not for any other set of variables.  Induct on $n$.
For  $n=2$, we use the previous lemma and Proposition~\ref{P_nPr} (a) to write
$$
\{p\}\{2\}_{s P_{2p}, t^p} = \{2p\} =P_2 P_p P_{2p}.
$$
Solving for $P_{2p}$ and using the fact that $\{p\}=P_p$ completes the base case.

For the induction step we use in turn Proposition~\ref{P_nPr} (a), the hypotheses on $p$ and $n$,  induction, and 
Lemma~\ref{pn:decomp} to obtain
\begin{align*}
\{pn\}
&= \prod_{d|pn} P_d\\[5pt]
&=\prod_{d|n} P_d P_{pd}\\[5pt]
&=P_p P_n P_{pn} \prod_{d|n \atop d\neq 1,n} P_d \cdot P_d(s P_{2p}, t^p)/P_d\\[5pt]
&=\frac{P_p P_n P_{pn} \{n\}_{s P_{2p}, t^p}}{P_n(s P_{2p}, t^p)}\\[5pt]
&=\frac{P_n P_{pn} \{pn\}}{P_n(s P_{2p}, t^p)}
\end{align*}
Solving for $P_{pn}$  finishes the proof.
\eprf

We can use this theorem to give a new relation between cyclotomic polynomials.  Note that  setting $s=q+1$ and $t=-q$ 
in the left-hand side of~\ree{ppm1} we get, using~\ree{lpq},
$$
\{p+1\}+t\{p-1\}= [p+1]_q-q[p-1]_q= q^p + 1.
$$
Using this substitution, we have the following immediate corollary of Theorem~\ref{pn}.
\bco
If $n\ge2$ is a positive integer and $p$ is prime not dividing $n$, then

\vs{10pt}
\eqed{
\Phi_{pn}(q) \Phi_n (q) = P_n(q^p+1,\ep(p)q^p).
}
\eco

We also have a Lucas analogue of Theorem~\ref{rf} (b).

\bth
\label{p^m n}
If $n$ is a positive integer, $p$ is a  prime not dividing $n$, and $m\ge2$ then
$$
P_{p^m n}(s,t) = \case{P_{p^{m-1} n}(s^2+2t, -t^2)}{if $p=2$,}{{P_{p^{m-1} n}(s P_{2p}, t^p)}}{if $p\ge3$.\rule{0pt}{15pt}}
$$
\eth
\bprf
We induct on $m$, where the base case is similar enough to the induction step that we will only provide  details for the latter.  And we will also just consider odd primes for similar reasons.  Given $m$, we induct on $n$.  For $n=1$, by Lemma~\ref{pn:decomp} we have
$$
\{p^m\}=\{p\}\{p^{m-1}\}_{s P_{2p}, t^p}
$$
Now expand both sides, using Proposition~\ref{P_nPr} (a) and use the fact that $\{p\}=P_p$, to get
$$
P_p P_{p^2} \cdots P_{p^m} 
= P_p \cdot P_p(s P_{2p}, t^p) \cdot P_{p^2}(s P_{2p}, t^p) \cdot\ \cdots\  \cdot P_{p^{m-1}}(s P_{2p}, t^p).
$$
Using the induction hypothesis on $m$ to cancel all but one factor on each side  gives the desired equality.
To deal with $n\ge2$, we do an induction on $n$ as well.
Expand $\{p^m n\}$ in a similar fashion  to what was done for $\{pn\}$ in Theorem~\ref{pn}.  After cancellation of terms, which uses the induction hypotheses on both $m$ and $n$, one obtains $P_{p^m n}/ P_{p^{m-1} n}(s P_{2p}, t^p)=1$ which is what we wish to prove. 
\eprf

Again, we can get a relation between cyclotomic polynomials and Lucas atoms by specialization.
\bco
If $n$ is a positive integer, $p$ is prime not dividing $n$, and $m\ge2$ then

\vs{10pt}
\eqed{
\Phi_{p^m n}(q)  = P_{p^{m-1} n}(q^p+1,\ep(p)q^p).
}
\eco

\section{Evalutations}
\label{e}

There are a number of interesting evaluations of the cyclotomic polynomials at various integers.  For example,
Suppose $b>1$ is an integer relatively prime to the prime $p$, and $n$ is the multiplicative order of $b$ modulo $p$.
Then it follows quickly from~\ree{PhiProd} that $p | \Phi_n(b)$.   For a more substantive example, there is the following conjecture which is implied by a conjecture of Bouniakowsky~\cite{bou:ntr}.
\bcon
For every positive integer $n$ there are infinitely many positive integers $b$ such that $\Phi_n(b)$ is prime.
\econ

We will prove some facts about the Lucas atoms modulo two and three.  The proofs will provide an application of the reduction formulas from Section~\ref{rfsec}.  They will also permit us to say something about the divisibility of the cyclotomic polynomials themselves.
We first need some information about the coefficients of $P_n(s,t)$.  
\ble
\label{P_n lead}
For $n\ge3$  we have
$$
P_n =\sum_{k=0}^{\phi(n)/2} c_k s^{\phi(n)-2k} t^k
$$
for certain constants $c_k$, where $c_0=1$ and 
$$
c_{\phi(n)/2}=\case{p}{if $n=2\cdot p^m$ for a prime $p\ge2$ and $m\ge1$,}{1}{else.}
$$
\ele
\bprf
All of these statements about $P_n$ are proved similarly, so we will just present a demonstration for the value of $c_{\phi(n)/2}$.
We induct on $n$, where the case $n=3$ is easy to check. 
From Lemma~\ref{cT(n)} we can write
\beq
\label{nExp}
\{n\}=\sum_{j=0}^{\fl{(n-1)/2}} a_j s^{n-j-1} t^j
\eeq
where the largest power of $t$ has coefficient
\beq
\label{tCoef}
a_{\fl{(n-1)/2}} =\case{n/2}{if $n$ is even,}{1}{if $n$ is odd.}
\eeq
 Now using Proposition~\ref{P_nPr} (a), induction, and the fact that $\sum_{d|n} \phi(d) = n$, we get from~\ree{nExp} that the  degree of $P_n$ as a polynomial in $t$ is $\phi(n)/2$.  
Using the same line of reasoning with~\ree{tCoef} we see that $c_{\phi(n)/2}=1$ for $n$ odd.
To complete the proof, we now repeat this argument in turn for the cases of $n=2\cdot p^m$ where $p$ is prime, and of $n = 2^l \cdot k$ where  $k$ is odd 
and either $l\ge2$ or $k$ has at least two prime factors.  The details are left to the reader. 
\eprf

We can now determine the behavior of $P_n(s,t)$ when $s,t$ are taken modulo $2$.

\bth
\label{PMd2}
Suppose $n\ge2$.  Then
\ben
\item[(a)] $P_n(0,0)=0$,
\item[(b)] $P_n(1,0)=1$,
\item[(c)] $2|P_n(0,1)$ if and only if $n=2^m$ for some $m\ge1$,
\item[(d)] $2|P_n(1,1)$ if and only if $n = 3\cdot 2^m$ for some $m\ge0$.
\een
\eth
\bprf 
The first three statements follow easily from the previous lemma.  So consider $P_n(1,1)$.
 Suppose first that $3$ does not divide $n$.
Let the $n$th Fibonacci number be denoted $F_n$ and recall that $F_n=\{n\}_{1,1}$.
It is well known and simple to prove that $2|F_n$ if and only if $3|n$.  So if $3$ is not a divisor of $n$ then $\{n\}_{1,1}$ is odd.  Thus the same must be true of its factor $P_n(1,1)$.

Since $P_3(1,1)=2$, we will now consider $n=3k$ where $k\ge2$ is not divisible by $3$.  
From Theorem~\ref{pn} we see that 
$P_{3k}(1,1)\Cong P_k(0,1)/P_k(1,1)\ (\Mod 2)$ since, as we have just proved, the denominator is not divisible by $2$.  By  part (c), 
$P_k(0,1)$ is even precisely when $k\ge2$ is a power of $2$, which finishes this case.

Finally, suppose $n=3^m k$ where $m\ge2$ and $k$ is not divisible by $3$.  By Theorem~\ref{p^m n} we have
$P_{3^m k}(1,1)\Cong P_{3^{m-1} k }(0,1)\ (\Mod 2)$.  But $3^{m-1}k$ is never a power of two since $m\ge2$.  So, by part (c) again, we have that $P_{3^{m-1} k }(0,1)$, and thus $P_{3^m k}(1,1)$, is odd as announced in the statement of the theorem.
\eprf

We can use the previous result to find the highest power of two which divides an evaluation of a cyclotomic polynomial.  For any prime $p$ and integer $n$ we let $\nu_p(n)$ be the highest power of $p$ dividing $n$.
\bco
\label{PhMd2}
If $b$ is an integer and $n\ge3$.  Then
$$
\nu_2(\Phi_n(b))=\case{1}{if $n=2^m$ for some $m\ge2$ and $b$ is odd,}{0}{else.}
$$
\eco
\bprf
We have $\Phi_n(b)=P_n(b+1,-b)$.   So we are only interested in the case where the two arguments in $P_n$ are of different parity.   But by Theorem~\ref{PMd2}, the only time $P_n(b+1,-b)$ for can be even for $n\ge3$ is when $n=2^m$ for some $m\ge2$.  So we need to investigate what happens when $\Phi_{2^m}(b)=b^{2^{m-1}}+1$.   Clearly if $b$ is even then this is not divisible by $2$.  And it is also easy to check that if $b$ is odd then, since $2^{m-1}$ is even, we have
$\Phi_{2^m}(b)\Cong 2\ (\Mod 4)$ which completes the proof.
\eprf

The proofs of the next two results are similar enough to those of Theorem~\ref{PMd2} and Corollary~\ref{PhMd2} that we will omit them.  However, as a labor-saving device, we note that because the powers of $s$ in $P_n(s,t)$ are all even for $n\ge3$, we always have $P_n(a,b)=P_n(-a,b)$.
\bth
\label{PMd3}
Suppose $n\ge3$.  Then
\ben
\item $P_n(0,0)=0$,
\item $P_n(\pm1,0)=1$,
\item $3|P_n(0,\pm1)$ if and only if $n=2\cdot 3^m$ for some $m\ge1$,
\item $3|P_n(\pm1,1)$ if and only if $n=4\cdot 3^m$ for some $m\ge0$,
\item $3|P_n(\pm1,-1)$ if and only if $n=3\cdot 3^m$ for some $m\ge0$.\hqed
\een
\eth

\bco
\label{PhMd3}
If $b$ is an integer and $n\ge3$.  Then

\vs{10pt}

$\rule{1ex}{0ex}\hfill{\dil \nu_3(\Phi_n(b))=\case{1}{if $n=3^m$ for some $m\ge1$ and $b\Cong 1\ (\Mod 3)$,}{0}{else.}}\hfill\raisebox{-10pt}{\qed}$
\eco

We note that, as opposed to the situation in Corollaries~\ref{PhMd2} and~\ref{PhMd3}, one can have $\nu_p(\Phi_n(b))\ge2$ for primes other than $2$ and $3$.  For example $\Phi_4(7)=50=2\cdot 5^2$.  We also remark that extending Theorems~\ref{PMd2} and~\ref{PMd3} to arbitrary primes is almost certainly hard.  One of the crucial tools in their proofs is the knowledge of the period of the Fibonacci sequence modulo $2$ and modulo $3$.  Although it is easy to see that this sequence is periodic modulo any integer, finding a formula for the period is a famous unsolved problem.

\section{Comments and open problems}

We will now present some avenues for future research hoping that the reader will be interested in exploring them.

\medskip

{\bf (1) Combinatorial interpretations.}  Since the Lucas atoms have nonnegative integer coefficients, one would hope that they count something.  But we have been unable to come up with a simple combinatorial interpretation for these polynomials, despite the fact that there are various well-known interpretations for the Lucas polynomials  themselves.  By using the reduction formulas, we have determined a complicated way of describing $P_n(s,t)$ when $n$ is a power of a prime in terms of certain colored tilings.  But it seems unlikely that this will extend to all $n$.  Once an interpretation is in place, it would be nice to take that as the {\em definition} of the Lucas atoms and then derive properties such as the decomposition~\ree{prod} combinatorially.

\medskip

{\bf (2) Alternating gamma vectors.}  One of the reasons for interest in gamma expansions is because of their connection with unimodality.  Call a polynomial $p(q)=\sum_{j\ge0}a_j q^j$ with real coefficients {\em unimodal} if 
$$
a_0\le a_1\le \ldots \le a_m\ge a_{m+1} \ge \ldots
$$
for some index $m$.  Unimodal sequences abound in algebra, combinatorics, and geometry.  See the survey articles of Stanley~\cite{sta:lus} and Brenti~\cite{bre:lus} and Br\"and\'en~\cite{bra:ulr} for more information.  Now suppose that $p(q)$ is palindromic.  If its gamma coefficients are all nonnegative, then $p(q)$ must be unimodal since all the polynomials involved in its expansion are unimodal with the same center of symmetry.
However, the definition of the map $\Ga$ in~\ree{GaDef} suggests that it might also be interesting to look at gamma expansions where the coefficients alternate in sign.  For example, this is true of the gamma expansions of the cyclotomic polynomials and their products.  Very little work has been done in this direction and we are only aware of a single paper of Brittenham, Carroll, Petersen, and Thomas~\cite{bcpt:uag} on this topic.

\medskip

{\bf (3) Coxeter groups.}  There are several ways in which the proofs of Theorems~\ref{CatThm} and~\ref{NarThm} could be improved.  First, it would be nice to have uniform proofs for all finite irreducible $W$ rather than having to go case-by-case.  It would also be desirable to find combinatorial proofs, especially in the cases where one is not already known.  And the best scenario would be to have these proofs rely on the combinatorics of the groups themselves.  In particular, it would be very interesting if these Lucas analogues are the generating functions for some statistics on the poset of noncrossing partitions $NC(W)$ which would reduce to the original counts when $s=2$ and $t=-1$.

\medskip

{\em Acknowledgement.}  We wish to thank Richard Stanley who originally had the idea of factoring the Lucas polynomials into Lucas atoms.  Without his seminal insight, this paper could not have been written.

\nocite{*}
\bibliographystyle{alpha}

\begin{thebibliography}{ABC{\etalchar{+}}15}

\bibitem[ABC{\etalchar{+}}15]{abcdl}
Farid Aliniaeifard, Nantel Bergeron, Cesar Ceballos, Tom Denton, and Shu~Xiao
  Li.
\newblock Algebraic {C}ombinatorics {S}eminar, {F}ields {I}nstitute.
\newblock \texttt{http://garsia.math.yorku.ca/fieldseminar/}, 2013--2015.

\bibitem[Arm09]{arm:gnp}
Drew Armstrong.
\newblock Generalized noncrossing partitions and combinatorics of {C}oxeter
  groups.
\newblock {\em Mem. Amer. Math. Soc.}, 202(949):x+159, 2009.

\bibitem[Ath]{ath:gcg}
Christos~A. Athanasiadis.
\newblock Gamma-positivity in combinatorics and geometry.
\newblock Preprint {\texttt{arXiv:1711.05983}}.

\bibitem[BCMS]{bcms:cil}
Curtis Bennet, Juan Carrillo, John Machacek, and Bruce Sagan.
\newblock Combinatorial interpretations of lucas analogues of binomial
  coefficients and catalan numbers.
\newblock Preprint {\texttt{arXiv:1809.09036}}.

\bibitem[BCPT16]{bcpt:uag}
Charles Brittenham, Andrew~T. Carroll, T.~Kyle Petersen, and Connor Thomas.
\newblock Unimodality via alternating gamma vectors.
\newblock {\em Electron. J. Combin.}, 23(2):Paper 2.40, 22, 2016.

\bibitem[Bou57]{bou:ntr}
Victor Bouniakowsky.
\newblock Nouveaux th\'eor\`emes relatifs \`a la distinction des nombres
  premiers et \`a la d\'ecomposition des entiers en facteurs.
\newblock {\em M\'em. Acad. Sc. St. P\'etersbourg}, 6:305--329, 1857.

\bibitem[BP09]{bp:caf}
Arthur~T. Benjamin and Sean~S. Plott.
\newblock A combinatorial approach to {F}ibonomial coefficients.
\newblock {\em Fibonacci Quart.}, 46/47(1):7--9, 2008/09.

\bibitem[Br\15]{bra:ulr}
Petter Br\"{a}nd\'{e}n.
\newblock Unimodality, log-concavity, real-rootedness and beyond.
\newblock In {\em Handbook of enumerative combinatorics}, Discrete Math. Appl.
  (Boca Raton), pages 437--483. CRC Press, Boca Raton, FL, 2015.

\bibitem[Bre94]{bre:lus}
Francesco Brenti.
\newblock Log-concave and unimodal sequences in algebra, combinatorics, and
  geometry: an update.
\newblock In {\em Jerusalem combinatorics '93}, volume 178 of {\em Contemp.
  Math.}, pages 71--89. Amer. Math. Soc., Providence, RI, 1994.

\bibitem[Ekh11]{ekh:ssl}
Shalosh~B. Ekhad.
\newblock The {S}agan-{S}avage {L}ucas-{C}atalan polynomials have positive
  coefficients.
\newblock Preprint {\texttt{arXiv:1101.4060}}, 2011.

\bibitem[Gau86]{gau:da}
Carl~Friedrich Gauss.
\newblock {\em Disquisitiones arithmeticae}.
\newblock Springer-Verlag, New York, 1986.
\newblock Translated and with a preface by Arthur A. Clarke, Revised by William
  C. Waterhouse, Cornelius Greither and A. W. Grootendorst and with a preface
  by Waterhouse.

\bibitem[Luc78a]{luc:tfn1}
Edouard Lucas.
\newblock Theorie des {F}onctions {N}umeriques {S}implement {P}eriodiques.
\newblock {\em Amer. J. Math.}, 1(2):184--196, 1878.

\bibitem[Luc78b]{luc:tfn3}
Edouard Lucas.
\newblock Theorie des {F}onctions {N}umeriques {S}implement {P}eriodiques.
\newblock {\em Amer. J. Math.}, 1(4):289--321, 1878.

\bibitem[Luc78c]{luc:tfn2}
Edouard Lucas.
\newblock Theorie des {F}onctions {N}umeriques {S}implement {P}eriodiques.
  [{C}ontinued].
\newblock {\em Amer. J. Math.}, 1(3):197--240, 1878.

\bibitem[Nen]{nen}
Gleb Nenashev.
\newblock personal communication.

\bibitem[Rie94]{rie:pnc}
Hans Riesel.
\newblock {\em Prime numbers and computer methods for factorization}, volume
  126 of {\em Progress in Mathematics}.
\newblock Birkh\"{a}user Boston, Inc., Boston, MA, second edition, 1994.

\bibitem[RS17]{rs:dsp}
Sujit Rao and Joe Suk.
\newblock Diheadral sieving phenomena.
\newblock Preprint {\texttt{arXiv:1710.06517}}, 2017.

\bibitem[SS10]{ss:cib}
Bruce~E. Sagan and Carla~D. Savage.
\newblock Combinatorial interpretations of binomial coefficient analogues
  related to {L}ucas sequences.
\newblock {\em Integers}, 10:A52, 697--703, 2010.

\bibitem[Sta89]{sta:lus}
Richard~P. Stanley.
\newblock Log-concave and unimodal sequences in algebra, combinatorics, and
  geometry.
\newblock In {\em Graph theory and its applications: East and West (Jinan,
  1986)}, volume 576 of {\em Ann. New York Acad. Sci.}, pages 500--535. New
  York Acad. Sci., New York, 1989.

\end{thebibliography}

\newcommand{\etalchar}[1]{$^{#1}$}

\end{document}